\documentclass{kybernetika}

\newtheorem{theorem}{Theorem}[section]
\newtheorem{lemma}[theorem]{Lemma}

\newtheorem{remark}[theorem]{Remark}

\def\deq{\stackrel{\Delta}{=}}    

\def\ra{\rightarrow}

\def\a{\alpha}

\def\g{\gamma}

\def\l{\lambda}
\def\L{\Lambda}

\def\O{\Omega}
\def\w{\infty}
\def\n{\nabla}

\def\R{\mathbb{R}}

\def\mc{\mathcal}
\def\mb{\mathbb}

\def\thefootnote{*}
\def\min{{\rm min}}

\def\be{\begin{array}}
\def\en{\end{array}}
\def\befc{ \begin{equation}\left\{\begin{array}{llll}}
\def\enfc{ \end{array} \right. \end{equation}}
\def\beeq{\begin{equation}}
\def\eneq{\end{equation}}
\def\beit{\begin{itemize}}
\def\enit{\end{itemize}}
\def\bepm{\begin{pmatrix}}
\def\enpm{\end{pmatrix}}

\def\thefootnote{$\dag$}

\def\hx{\mbox{\boldmath$x$}}   
\def\hy{\mbox{\boldmath$y$}}   

\def\h1{\mbox{\boldmath$\mathbbm 1$}}

\def\hlbb{\mbox{\boldmath$\lbb$}}   

\DeclareMathOperator{\argmin}{argmin}
\DeclareMathOperator{\prox}{Prox}
\usepackage{bbm}
\usepackage[T1]{fontenc}
\usepackage{mathrsfs}
\usepackage{amsmath,multicol,amsopn}
\usepackage{cases}
\usepackage{latexsym,amssymb}
\usepackage{amsbsy, amstext}
\usepackage{graphicx}    
\usepackage{empheq}
\newtheorem{assumption}{Assumption}
\allowdisplaybreaks[4]   
\numberwithin{equation}{section}

\hypersetup{
  unicode=true,
  colorlinks=true,
  linkcolor=blue,
  citecolor=blue,
  urlcolor=blue
}

\makeatletter
\newcommand{\lb}{{\mathchoice
  {\smash@bar\textfont\displaystyle{0.45}{1.2}\lambda}
  {\smash@bar\textfont\textstyle{0.25}{1.2}\lambda}
  {\smash@bar\scriptfont\scriptstyle{0.25}{1.2}\lambda}
  {\smash@bar\scriptscriptfont\scriptscriptstyle{0.25}{1.2}\lambda}
}}
\newcommand{\Lb}{{\mathchoice
  {\smash@bar\textfont\displaystyle{0.45}{1.2}L}
  {\smash@bar\textfont\textstyle{0.25}{1.2}L}
  {\smash@bar\scriptfont\scriptstyle{0.25}{1.2}L}
  {\smash@bar\scriptscriptfont\scriptscriptstyle{0.25}{1.2}L}
}}
\newcommand{\lbb}{{\mathchoice
  {\smash@bar\textfont\displaystyle{0.15}{1.2}\smash@bar\textfont\displaystyle{0.40}{1.2}\lambda}
  {\smash@bar\textfont\textstyle{0.15}{1.2}\smash@bar\textfont\displaystyle{0.45}{1.2}\lambda}
  {\smash@bar\scriptfont\scriptstyle{0.15}{1.0}\smash@bar\textfont\displaystyle{0.75}{0.5}\lambda}
  {\smash@bar\scriptscriptfont\scriptscriptstyle{0.15}{1.0}\smash@bar\textfont\displaystyle{0.75}{0.5}\lambda}
}}

\newcommand{\smash@bar}[4]{%
  \smash{\rlap{\raisebox{-#3\fontdimen5#10}{$\m@th#2\mkern#4mu\mathchar'26$}}}%
}
\makeatother

\newenvironment{assbis}
  {\addtocounter{assumption}{-3}%
   \begin{assumption}}
  {\end{assumption}}

\renewcommand{\l}{\lambda}
\renewcommand{\L}{\Lambda}

\begin{document}
\pagestyle{myheadings}

\title{A projection-free dynamics for nonsmooth\\ composite optimization}

\author{Wei Ni, Yangfan Qiu, and Yanyan Xiao}

\contact{Wei}{Ni}{School of Mathematics and Computer Sciences, Nanchang University, Nanchang, China. He is also with Key Laboratory of Intelligent Systems and Human-Machine Interaction of Jiangxi Province.}{niwei@ncu.edu.cn}
\contact{Yanyan}{Xiao}{School of Mathematics and Computer Sciences, Nanchang University, Nanchang, China. }{15970167278@163.com}
\contact{Yangfan}{Qiu}{School of Mathematics and Computer Sciences, Nanchang University, Nanchang, China.}{875877601@qq.com}

\markboth{W. Ni et al. } {A projection-free dynamics for nonsmooth composite optimization}

\maketitle

\begin{abstract}
This paper proposes a projection-free primal-dual dynamics for the nonsmooth composite optimization problems with equality and inequality constraints.
To deal with optimization constraints, this paper departs from the use of gradient projection method,  but resorts to the idea of mirror descent to design a continuous-time smooth optimization dynamics which advantageously leads to easier convergence analysis and more efficient numerical simulation.
Also, the strategy of proximal augmented Lagrangian (PAL$^{\dag}$) is extended to incorporate general convex equality-inequality constraints and the  strong convexity-concavity of the primal-dual variables is achieved, ensuring exponential convergence of the resulting algorithm.
Furthermore, the convergence result in this paper extends existing exponential convergence  which  either takes no account of  constraints or considers only affine linear constraints, and it also enhances  existing asymptotic convergence under convex constraints which unfortunately  depends on the complex gradient projection scheme.
\end{abstract}



\let\thefootnote\dag\footnotetext
{For convenience, other acronyms defined in this paper are also summarized here:
\vspace{-0.3cm}
\begin{align*}
\begin{array}{lllllllll}
&\text{AL: augmented Lagrangian.} &  \text{PAL: proximal augmented Lagrangian.}\\
&\text{IQC: integral quadratic constraints. }  & \text{SCQ: Slater's constrained qualification. } \\
& \text{KKT: Karush-Kuhn-Tucker. } &    \text{LICQ: linear independence constraint qualification. }
\end{array}
\end{align*}
}

\section{INTRODUCTION}
Composite optimization with nonsmooth regularization has been a topic of a large amount of research. This problem consists in minimizing the sum of a differentiable function $f(x)$ and a possibly nonsmooth one $\varphi(Tx)$, with the latter being a composition of an extended-valued function $\varphi$ with a linear operator $T$.   A prominent example is the empirical risk minimization with nonsmooth regularization as parameter choice rule in deep neural network \cite{lecun2015}.  As another example, the typical signal processing problem of  restoration of a continuous signal from its digital counterparts is tackled by solving an optimization problem known as the  lasso  \cite{tibshirani1996} which is a specific case of the problem considered in this paper. Applications in other fields include inverse problem \cite{beck2009},  statistical estimation \cite{bickel2008}, and control \cite{lin2013}.  Due to the presence of the nondifferentiable part $\varphi$ in the objective function, classical gradient based optimization algorithms,  such as gradient descent, conjugate gradient, and Newton methods, cannot be directly applied. To tackle non-smoothness, a direct  research line  is the nonsmooth optimization method, including subgradient method \cite{shor2012}, stochastic approximation \cite{nemirovski2009} and bundle method \cite{makela2002}, but the convergence is very slow.
 As a parallel, the smooth optimization method \cite{nesterov2005} replaces the nonsmooth component with smooth approximations parameterized by a positive parameter and then resorts to smooth optimization routines; however, as the parameter approaches zero, the problem becomes increasingly involved.
 A breakthrough was made via proximal operator in a series of milestone papers  which do not require explicit gradient \cite{rockafellar1976b,beck2009,boyd2011,tseng2000}. In these methods, the smooth part of the objective function is processed by evaluation of its gradient operator, while a nonsmooth part is handled via its Moreau proximal operator \cite{moreau1965}, resulting in iterations of the form $x^{+} \in \operatorname{prox}_{\gamma \psi (T(x))}(x-\gamma \nabla f(x))$  for some suitable stepsize $\g$.
 Proximal method can be derived from operator splitting techniques (see Douglas \cite{douglas1956}, Passty \cite{passty1979} and Bruck \cite{bruck1977}) which find their roots in the works of Browder \cite{browder1965} and Minty \cite{minty1962},  and now it is developed in the more general setting of maximal monotone operators \cite{abbas2015}.  The more recent work by Combettes and Wajs \cite{combettes2005} provides important insights into the proximal method and has popularized the method for a wide audience.  Nowdays the proximal method becomes an important part of convex analysis, optimization, and nonlinear functional analysis.

 The above proximal method involves the proximal calculation with respect to the composite function $\varphi\circ T$ which is however difficult to evaluate for general case.  To overcome this, a frequently used technique is the variable splitting which constructs a new variable $y$ to serve as the argument of $\varphi$ and imposes additional constraint $y=Tx$.
Then this constrained optimization problem can be tackled by using an augmented Lagrangian (AL) method \cite{nocedal1999} originally due to Hestenes \cite{hestenes1969} and Powell \cite{powell1969} and by designing corresponding primal-dual algorithms \cite{arrow1958}.
The AL here extends the classical Lagrangian by including an extra quadratic term,  without changing the optimal solution but facilitating the proximal operator to be used here.  The proximal version of this method is earlier  developed by  Rockafellar  \cite{rockafellar1976b} and its connection to  alternating direction of multiplier method is recently explored by Boyd \cite{boyd2011},  followed by  a burst of research activities in subsequent years.
Note that the AL function mentioned above still contains a nondifferentiable term and thus precludes the standard strategy of gradient descent on  primal variables and gradient ascent on dual variables. To conquer the nondifferentiablity in the traditional AL, Dhingra et al. \cite{dhingra2019} recently propose a new kind of differentiable Lagrangian, called PAL,  which is obtained by taking infimum over the AL with respect to the primal variable of the nonsmooth function. Compared with the AL in \cite{nocedal1999}, a coming advantage of the PAL  is that it is continuously differentiable on both primal and dual variables so that direct primal-dual dynamics can be applied on the PAL.

With the help of PAL, the original optimization problem is transformed into the computation of the saddle point of this PAL,  and various formulations of primal-dual dynamics have been explored in the literature with focus on investigating asymptotic  or exponential convergence under different conditions \cite{dhingra2021,goldsztajn2021,hassan2021,ju2022,francca2021,tang2020,eckstein2025}. In a broader view, ever since von Neumann's celebrated min-max theorem \cite{neumann1928}, the problem of finding a saddle point for a convex-concave function has been a major focus of disciplines like machine learning \cite{goodfellow2014}, robust optimization \cite{ben2009}, robust control \cite{hast2013} and game theory \cite{bacsar1998}. For the abstract gradient descent-ascent algorithm which  performs gradient descent in convex variable and gradient ascent in concave variable,   there is a vast literature obtaining  asymptotic convergence  if either the convexity or concavity  are strict \cite{sion1958,kose1956,arrow1958,cherukuri2017}; exponential convergence can be further ensured if strongly convex strongly concave condition is imposed \cite{cherukuri2017}.
However, when specialized to the case of classical Lagrangian which is only linear  on the dual variables, exponential stability can not  be obtained directly for the primal-dual dynamics. This has also been highlighted in \cite{wang2021} that  the combination of general convex constraints and nonsmooth objectives poses a significant challenge to achieving exponential stability.
To prove exponential stability under extra conditions, other routines  include the contraction analysis \cite{cisneros2021} and the technique of  spectral bounds of saddle matrices \cite{benzi2005}. In this paper, instead of building on the classical Lagrangian but on the PAL, we develop a primal-dual optimization dynamics.  Aside from the advantage of variable splitting inherent in the PAL,  we will show that the PAL induces strong concavity on dual variable under appropriate conditions and renders the corresponding primal-dual optimization dynamics to achieve exponential stability by taking the method of  Lyapunov function and convex analysis. Obviously, improved upon asymptotic stability, exponential stability offers explicit convergence rate of the optimization algorithms.

When considering constraints in the PAL method, most existing results have been successful in proving exponential stability either for equality constraints \cite{wang2021,cisneros2021} or for (affine) linear inequality constraints \cite{adegbege2021,qu2019}.
When only asymptotic stability is required, general convex inequality constraints can be tackled  either by making  use of the framework of hybrid dynamical system \cite{feijer2010} or by referring to techniques involving projection \cite{goldsztajn2021,tang2020}, where the former  invokes a generalized Lasalle invariance principle developed for hybrid systems and the latter  utilizes the sophisticated nonsmooth analysis.
This generalized Lasalle invariance principle is also improved to be used for the projection approach; see Cherukuri et al. \cite{cherukuri2016},  \cite{cherukuri2018}. Unfortunately, both the hybrid system approach and the projection approach  would result in discontinuity of the primal-dual dynamics and could bring difficulties in the convergence proof.
To overcome this difficulty, this paper pursues a projection free method to design a smooth primal-dual optimization dynamics  by utilizing the technique of mirror descent \cite{nemirovsky1983,raginsky2012} under the general framework of barrier operator \cite{bolte2003}. Our
method avoids projection and thus reduces the difficulties of convergence analysis and
iterative computation.

Aside from merit of yielding a smooth optimization dynamics, further advantage of our method lies in its ability to deal with general nonlinear convex optimization constraints.
For optimization algorithm design and corresponding exponential convergence analysis, a popular choice is the IQC framework \cite{lessard2016,fazlyab2018} which provides  linear matrix inequality  based conditions for exponential stability of  the feedback interconnection of a linear system with a nonlinearity satisfying IQC. Although the gradient mapping or the proximal mapping individually satisfies IQC, composition of each with a nonlinear map generally does not satisfy IQC; instead, IQC holds for compositions with affine operations \cite{fazlyab2018}. In view of this fact, most exponential convergence of PAL methods  based on IQC deals with either equality constraints  \cite{dhingra2019,hassan2020,hassan2021,wang2021,ding2020} (which should be affine in decision variables in view of the convex nature of the problem) or affine inequality constraints \cite{hassan2018,bansode2019}, with general nonlinear convex inequality constraints not touched.  Similar situations occur in other literature;  for one example  in \cite{adegbege2021}, the authors proposed another projection-free method to deal with inequality constraints, but only limited to affine inequality constraints and for another in \cite{qu2019}, exponential stability is obtained but limited only to affine inequality constraints. Our method does not rely on IQC analysis and exponential stability can be achieved under general convex inequality constraints.
Although general convex inequality constraints are considered in \cite{goldsztajn2021} and \cite{tang2020,bin2024}, the former can only guarantees asymptotic convergence and the later attains  semi-globally exponential convergence. These deficiencies of asymptotic convergence and semi-global exponential convergence are improved to exponential convergence in this paper.
We also note  the works \cite{qu2019,adegbege2021,bin2024} which handle inequality constraints without explicit projection. However,  \cite{qu2019,adegbege2021,bin2024} replace projection-based updates with max-based regularization (specifically algorithm (9a)-(9c) in \cite{adegbege2021}, (9a)-(9b) in \cite{qu2019}, and (1b) in \cite{bin2024}), implicitly ensuring nonnegativity of the multiplier associated with inequality constraint. While this improves computational simplicity, the max operator still makes the optimization dynamics nonsmooth in the classical sense.
In contrast, our  projection-free dynamics  is fully smooth, and this smooth nature enables rigorous stability certification via classical Lyapunov techniques, bypassing the need for nonsmooth analysis tools like differential inclusions.
More important, the max-based method in \cite{qu2019,adegbege2021}  are limited to affine inequality constraints, with no extension to general convex inequality constraints. Our approach overcomes this limitation and  extends to general convex inequality constraints.

To summarize, this paper proposes a projection-free method for nonsmooth composite optimization with constraints. The contributions of this paper are  as follows.
\vspace{-0.1cm}
\begin{itemize} \setlength{\parskip}{-4pt}
\item Firstly, equality and general inequality optimization constraints are taken into the nonsmooth composite optimization problem and a projection-free optimization dynamics is proposed. Unlike existing  projection or  max based methods of dealing with constraints which render nonsmooth optimization dynamics, the  optimization dynamics proposed in this paper is smooth,  with the benefit of reducing the difficulties of convergence analysis and iteration complexity,  bypassing the need for nonsmooth analysis tools like differential inclusions.
\item Secondly, improved on the classical Lagrangian which is only linear on the dual variable, the PAL method  is adopted and it  extends existing PAL method by considering optimization constraints.
    Strong concavity on dual variable is obtained in this paper so that exponential stability is achieved, thereby improving upon existing results that offer only asymptotic stability. In contrast  to asymptotic stability, exponential stability provides explicit convergence rates for the optimization algorithms.

\item Lastly, the method of achieving exponential stability in this paper differs existing IQC-based approaches that work only for  equality constraints or affine inequality constraints. our method does not rely on IQC analysis and is capable of handling general nonlinear convex inequality constraints while still ensuring exponential convergence,   going beyond the limitations of existing IQC-based approaches.
\end{itemize}

The rest of this paper is structured as follows. Preliminaries are put in Section \ref{prel}.  Section \ref{pf} describes the nonsmooth and composite optimization problem. Section  \ref{odca} is devoted the design of the optimization dynamics and the analysis of algorithm convergence. A simulation example is introduced in Section \ref{example} and a brief conclusion is made in Section \ref{conclusion}.

\section{Preliminaries} \label{prel}

We use $\mb R^n$ to denote the set of $n$-dimensional real-valued vectors and $\R^n_{+}$ those vectors with nonnegative components. For $a\in \R^n$, by $a \prec 0$ ($a\preccurlyeq 0$) we mean that each entry of $a$ is less than (less than or equal to) zero. Euclidian norm of a vector is denoted by $\|\cdot\|$.
Letting \scalebox{0.9}[1]{$a=(a_1, \cdots, a_n)^{\top}\in \mathbb{R}^n$} and \scalebox{0.9}[1]{$b=(b_1, \cdots, b_n)^{\top}\in \mathbb{R}^n$}, we define $a\odot\ b=(a_1b_1, \cdots, a_nb_n)^{\top}$ and $a\oslash b\!=\!(a_1/b_1, \cdots, a_n/b_n)^{\top}$. For vectors $\a_1, \cdots, \a_m$, the notation ${\rm col}\{\a_i\}_{i=1}^m$ denotes a new vector $(\a_1^{\top}, \cdots, \a_m^{\top})^{\top}$.

An extended value function $\varphi: \mathbb{R}^{m} \rightarrow \mathbb{R} \cup\{+\infty\}$ with effective domain $\operatorname{dom} (\varphi) \equiv\{y \mid \varphi(y)<\infty\}$ is proper if $\varphi(y)<+\infty$ for at least one $y\in \operatorname{dom}(\varphi)$ and $\varphi(y)>-\w$ for all $y\in \operatorname{dom}(\varphi)$.
A function $\varphi: \R^m \ra [-\w,\w]$ is lower semicontinuous if $\{\varphi\leq c\}\deq \{y\in \R^m: \varphi(y)\leq c\}$ is a closed subset of $\R^m$ for all $c\in \R$.
A function $\varphi: \mathbb{R}^{m} \rightarrow \mathbb{R} \cup\{+\infty\}$ is closed if its epigraph $\left\{(y, c) \in \mathbb{R}^{m} \times \mathbb{R} \mid \varphi(y) \leq c\right\}$ is a
closed set. For a convex and lower semicontinous function $\varphi$, its conjugate function (or Fenchel-Legendre transform)  $\varphi^*: \mathbb{R} \rightarrow \mathbb{R} \cup\{+\infty\}$ is defined as $\varphi^*(v)=\sup _{u \in \operatorname{dom}(\varphi)}\{u v-\varphi(u)\}$, which is also convex and lower-semicontinuous. Clearly only those $u$ in $\operatorname{dom}(\varphi)$ are relevant in the calculation of this supremum.

For an extended value function $\varphi: \mathbb{R}^{m} \rightarrow \mathbb{R} \cup\{+\infty\}$, its
subdifferential at the point $y\in \operatorname{dom}(\varphi)$  is denoted  as $\partial \varphi(y)=\{v \in \mathbb{R}^{y} \mid \varphi(z)-\varphi(y) \geq\langle v, z-y\rangle \text { for all } z\}$
which is always closed and convex for $y \in \operatorname{int}(\operatorname{dom} \varphi)$;  if $\varphi$ is convex then it is nonempty and if $\varphi$ is differentiable then $\partial \varphi(y)=\{\n \varphi(y)\}$, the usual gradient.

A differentiable function $f$ is $m_f$-strongly convex   if for any $x$ and $\hat x$, $f(\hat{x}) \geq f(x)+\langle\nabla f(x), \hat{x}-x\rangle+\frac{m_{f}}{2}\|\hat{x}-x\|_{2}^{2}$,
or equivalently, $\|\nabla f(\hat{x})-\nabla f(\hat{x})\|_{2} \geq m_{f}\|x-\hat{x}\|_{2}$,
or equivalently,$x\ra f(x)-\frac{m_f}{2}\|x\|^2$ is convex.
The gradient of a continuously-differentiable function $f$ is
Lipschitz continuous with parameter $\ell_f$ if for any $x$ and $\hat x$, $f(\hat{x}) \leq f(x)+\langle\nabla f(x), \hat{x}-x\rangle+\frac{\ell_{f}}{2}\|\hat{x}-x\|_{2}^{2}$,
or equivalently, $\|\nabla f(x)-\nabla f(\hat{x})\|_{2} \leq \ell_{f}\|x-\hat{x}\|_{2}$.
Moreover, if an $m_f$-strongly convex function $f$ has an
$\ell_f$-Lipschitz continuous gradient, then any $x$ and $\hat x$ it holds $\langle\nabla f(x)-\nabla f(\hat{x}), x-\hat{x}\rangle \geq \frac{m_{f} \ell_{f}}{m_{f}+\ell_{f}}\|x-\hat{x}\|_{2}^{2}+\frac{1}{m_{f}+\ell_{f}}\|\nabla f(x)-\nabla f(\hat{x})\|_{2}^{2}$.
We say that $f$ is strongly concave if $-f$ is strongly convex.

Let $K \subseteq \mathbb{R}^{n+m}$  be nonempty, closed and convex set. We say that a function $\varphi(x, y): K \rightarrow \mathbb{R}$ is concave-convex on $K$ if for any $\left(x, y\right) \in K$, $\varphi\left(x, y\right)$ is a concave function of $y$ and $\varphi(x, y)$ is a convex
function of $x$. If either the concavity or convexity is always
strict, we say that  is strictly concave-convex on $K$.
For a concave-convex function $\varphi: \mathbb{R}^{n} \times \mathbb{R}^{m} \rightarrow \mathbb{R}$ , we say that $(\bar{x}, \bar{y}) \in \mathbb{R}^{n+m}$  is a saddle point
if for all $x \in \mathbb{R}^{n} \text { and } y \in \mathbb{R}^{m}$ we have the inequality $\varphi(x, \bar{y}) \leq \varphi(\bar{x}, \bar{y}) \leq \varphi(\bar{x}, y)$.  If $\varphi$  is in addition continuously differentiable, then $(\bar x, \bar y)$ is a saddle point if and only
if $\n_{\bar x} \varphi(\bar{x}, \bar{y})=0 \text { and } \n_{\bar y}\varphi(\bar{x}, \bar{y})=0$.

For a non-differentiable extended value function $\varphi: \mathbb{R}^{m} \rightarrow \mathbb{R} \cup\{+\infty\}$,  its Moreau envelope $\varphi_{\mu}$ is defined as
\begin{align}\label{mor}
\varphi_{\mu}(v)\deq \underset{y\in\R^m}{\min }\left\{\varphi(y)+\frac{1}{2\mu}\|y-v\|^2 \right\} \leq \varphi(v).
\end{align}
Even though $\varphi$ may not be differentiable, its Moreau envelope may be, with $\operatorname{dom} (\varphi_{\mu})=\R^m$. We also note that $\varphi_{\mu}$ is convex since the function in the curly brackets is jointly convex in $(x, y)$ and thus the epigraph  of $\varphi_{\mu}$ is the projection of a convex set. The Moreau envelope $\varphi_{\mu}$ provides a smooth approximation of $\varphi$ from below and retains the same minimizers  \cite{rockafellar1976}. The minimizer of the  problem \eqref{mor}
is called as the proximal operator which is defined and calculated as follows,
\begin{align}
\prox_{\mu\varphi}(v)
&=\underset{y\in\R^m}{\argmin }\left\{\varphi(y)+\frac{1}{2\mu}\|y-v\|^2 \right\}\label{p1}\\
& =y-\mu \n \varphi_{\mu}(y), \label{p2}
\end{align}
where the second equality can be seen by noting that the gradient of the Moreau is
\begin{align*}
\n \varphi_{\mu}(y)=\frac{1}{\mu}[y-\prox_{\mu \varphi}(v)].
\end{align*}
The proximal operator expressed in the equality \eqref{p2}  can be viewed as a gradient descent for the Moreau envelope $\varphi_{\mu}$.   It can  also be viewed as a  generalization of  the notion of Euclidean projection if one takes $\varphi$ to be the indicator function of some convex set $C$ so that $\prox_{\mu \varphi}(y)=\argmin_{z\in C}\|y-z\|^2$.
For some choices of $\varphi$, the corresponding $\prox_{\mu\varphi}$ has  well known closed forms. For example, if $\varphi(x)=\|x\|_{1}$, then $\prox_{\mu\varphi}=$ $\operatorname{soft}(\mathrm{y}, \mu)$, where $\operatorname{soft}(\cdot, \mu)$ denotes the component-wise application of the soft-threshold function $y \mapsto \operatorname{sign}(y) \max \{|y|-$ $\mu, 0\}$.
In general the proximal operator may be set-valued, but it is scalar-valued
if $f(x)$ is  proper and convex.

\section{Problem formulation}\label{pf}

Consider the following nonsmooth composite optimization problem with equality and inequality constraints:
\begin{subequations} \label{optimization}
    \begin{empheq}[left={\mathcal{P}: \empheqlbrace\,}]{align}
  &\underset{x\in\R^n}{\rm minimize}  &    f(x)+\varphi(Tx), \label{objective}\\
 & {\rm subject~to}                              &     g(x)\preccurlyeq 0,  \label{c1}\\
 & {\rm }                                                &     h(x)=0, &&\label{c2}
    \end{empheq}
\end{subequations}
where $x\in \R^n$ is the decision variable, $T\in \R^{m\times n}$ is a given matrix of full column rank,   $f: \mathbb{R}^n \rightarrow \mathbb{R}$ and $\varphi: \mathbb{R}^m \rightarrow \mathbb{R} \cup \{+\w\}$ are proper, convex, and lower semicontinuous, with $f$ continuously differentiable  and $\varphi$ non-differentiable,
$g=(g_{1}, \cdots, g_{r})^{\top}: \mathbb{R}^n \rightarrow \mathbb{R}^{r}$,  and
$h=(h_{1}, \cdots, h_{s})^{\top}:\mathbb{R}^n \rightarrow \mathbb{R}^{s}$  are respectively the  inequality and the equality constraints, with  $g_j, j=1, \cdots, r$  being convex and continuously differentiable and $h_j, j=1, \cdots, s$ being affine. This makes the optimization problem \eqref{optimization} a convex one. In the problem \eqref{optimization}, the $T$  can represent a convolution operator for the image deconvolution problem  and represents operation of tomographic projections in the general image reconstruction problems \cite{afonso2010}.

By introducing an auxiliary variable $y\deq Tx$ and by imposing additional constraint $y=Tx$,  the original optimization problem \eqref{optimization} can be equivalently transformed  into the following one
\begin{subequations} \label{optimization2}
    \begin{empheq}[left={\widetilde{\mathcal{P}}: \empheqlbrace\,}]{align}
&\underset{x\in\R^n, y\in\R^m}{\rm minimize}   & \hspace{-4cm}  f(x)+\varphi(y), \label{objective}\\
& {\rm subject~to}                                     &   g(x)\preccurlyeq 0, \label{cc1}\\
&  {\rm }                                                      & h(x)=0, \label{cc2}\\
&                                                                  & Tx=y. \label{cc3}
    \end{empheq}
\end{subequations}
Clearly, $(x^*, y^*)$ is an optimal solution to \eqref{optimization2} if and only if $x^*$ is an optimal solution to the original optimization  \eqref{optimization}. The above procedure of creating a new variable $y$ to serve as the argument of $\varphi$ and imposing the constraint $y=Tx$ is a standard technique called variable splitting.  As will be seen, the transformed optimization problem \eqref{optimization2} is easier to solve than the original one \eqref{optimization}. The following assumptions are made throughout this paper.

\begin{assumption}\label{sconvex}
The function $f$ is $\a$-strongly convex for some $\a>0$ and continuously differentiable.
\end{assumption}

\begin{assumption}\label{lcontinuous}
The function $\varphi$ is a convex and subdifferentiable function whose subgradient is  $1/\ell$-Lipschtz continous for some $\ell>0$; that is,  for all $y, \tilde y \in \R^m$, it has $\left[d(y)-d(\tilde y)\right]^{T}(y-\tilde y) \leq 1/\ell \|y-\tilde y\|^{2}$, where $d(y)\in \partial g(y)$ and $d(\tilde y)\in \partial g(\tilde y)$.
\end{assumption}

Let $(x^*, y^*)$ be an optimal solution to the problem \eqref{optimization2}.  By assuming  additional assumptions on the constraint functions,  called constrained qualifications,   the following Karush-Kuhn-Tucker (KKT) systems hold at the minimizer  $(x^*, y^*)$: there exist multipliers
$(\lambda^*, \lb^*, \lbb^*)\in \mathbb R^{r}_{+} \times \R^s \times \R^m$  such that
\begin{subequations} \label{KKT2}
    \begin{empheq}[left={KKT: \empheqlbrace\,}]{align}
      & g(x^*) \preceq 0,  \label{KKT2a}\\
      & h(x^*) = 0, \label{KKT2b}\\
      & Tx^*=y^*, \label{KKT2c}\\
      & \lambda^* \succeq 0,  \label{KKT2d}\\
      &  \lambda^* \odot g(x^*)=0, \label{KKT2e}\\
      &\n f(x^*)\!+\!\partial \phi(y^*)\!+\!\sum\nolimits_{i=1}^r \l^{*}_{i} \n g_i(x^*) +\!\sum\nolimits_{i=1}^s \lb^{*}_{i} \n h_i(x^*)\!+\!T^{\top}\lbb^{*}\ni 0. \label{KKT2f}
    \end{empheq}
\end{subequations}
A widely used constrained qualification is the Slater's constrained qualification (SCQ):  there exist $(x, y)\in \mathbb{R}^n \times \R^m$ such that $ g(x) \!\prec \!0$, $h(x)\!=\!0$ and $Tx\!=\!y$.  In other  words, assuming SCQ,
``$(x^*, y^*)$ solves ($\widetilde{\mc P}$)'' $\Rightarrow$ ``$\exists$ a set of $(\lambda^*, \lb^*, \lbb^*)$ together with $(x^*, y^*)$ solving ($\mc K\mc K\mc T$)''. Furthermore, for convex problem, this implication is bidirectional. Refer to   \cite{boyd2004} for details.  Although SCQ guarantees existence of multipliers satisfying the  $\mc K\mc K\mc T$ system \eqref{KKT2}, it does not ensure uniqueness.
A stronger constrained qualification, named the linear independence constraint
qualification (LICQ), can achieve this goal.
Using
$\n_{\hspace{-0.1cm}J} g(x^*)$ to denote the submatrix of
$\n g(x^*)$ given by rows with indices in
$J(x^*)\!=\! \left\{i | g_i\left(x^{*}\right)\!=\!0\right\}$,
the LICQ is defined  as
\begin{align}\label{LICQ}\renewcommand\arraystretch{0.8}
\operatorname{rank}\left[
\begin{array}{lll}
\nabla h(x^*) \\
\nabla_{\hspace{-0.1cm}J}g(x^*)
\end{array}\right]=s+|J(x^*)|,
\end{align}
here $|\cdot|$ denotes the set cardinality.
By assuming LICQ,  one obtains a result \cite{wachsmuth2013} on the existence and uniqueness of multipliers satisfying \eqref{KKT2},
\begin{align}\label{uniqueMultiplier}
\scalebox{0.85}[1]{``$(x^*, y^*)$ solves ($\widetilde{\mc P}$)''} \Rightarrow \scalebox{0.85}[1]{``$\exists!$  $(\lambda_{ij}^*, \lb_{ij}^*, \lbb_{ij}^*)$ together with $(x^*, y^*)$ solving (KKT)''},
\end{align}
where the notation $\exists!$ stands for "there exists a unique".
This direction is bidirectional if the problem is convex. In view of above considerations, the following assumption is made.

\begin{assumption}\label{licq}
The constraint  functions in \eqref{cc1}-\eqref{cc2} satisfy LICQ.
\end{assumption}

\section{Optimization dynamics and convergence analysis} \label{odca}

This section proposes an exponentially stable primal-dual dynamical system without using gradient projection to solve the optimization problem \eqref{optimization} or \eqref{optimization2} building on the construction of a PAL.

\subsection{Design of PAL}
For the optimization problem \eqref{optimization2}, define the following AL function
\begin{align}\label{Lag1}
\hspace{-0.3cm} L^{\mu}\!(x, y; \l, \lb, \lbb)
\!=\!f(x)\!+\!\varphi(y)\!+\!\l^{\!\top}\! g(x)\!+\!\lb^{\!\top}\!h(x)\!+\!\lbb^{\!\top}\!(Tx\!-\!y)\!+\!\frac{1}{2\mu}\!\|Tx\!-\!y\|^2,
\end{align}
where $(x,y)\in \R^{n} \times \R^m$ are the primal variables and $(\l, \lb, \lbb)\in \R^r_{+} \times \R^s \times \R^m$ are the dual variables.
 Under the LIQC conditions in Assumption \ref{licq}, $(x^*, y^*)$ is an optimal solution to the optimization problem \eqref{optimization2} if and only if  there exists $(\l^*, \lb^*, \lbb^*)\in \R^r_{+} \times \R^s \times \R^m$  such that the following saddle point inequality
\begin{align}\label{saddle1}
 L^{\mu}(x^*, y^*; \l, \lb, \lbb)\leq L^{\mu}(x^*, y^*; \l^*, \lb^*, \lbb^*)\leq   L^{\mu}(x, y; \l^*, \lb^*, \lbb^*)
\end{align}
holds for all $(x, y)\in \R^n\times \R^m$ and $(\l, \lb,\lbb)\in \R^r_{+}\times \R^s \times \R^m$.

The saddle point property \eqref{saddle1}  leads us naturally to consider designing naive primal-dual dynamics based on the augmented Lagrangian $ L^{\mu}$ with gradient descent on $(x, y)$ and gradient ascent on $(\l, \lb, \lbb)$, with the aim to seek the optimal solution $(x^*, y^*)$ of the optimal problem \eqref{optimization2}. However, this naive primal-dual dynamics can only achieve asymptotic stability since $L^{\mu}$ is only linear on $(\l, \lb, \lbb)$; to achieve exponential stability of primal-dual dynamics, strong concavity  on $(\l, \lb, \lbb)$ should be required. Toward this goal, the subsequent paragraph introduces a PAL which is strongly concave in $(\l, \lb, \lbb)$.

Note that the AL in \eqref{Lag1} can be rewritten via completing of squares  as
\begin{align*}
 L^{\mu}(x, y; \l, \lb, \lbb)
\!=\!f(x)\!+\!\varphi(y)\!+\!\l^{\!\top}\! g(x)\!+\!\lb^{\!\top}\!h(x)\!+\!\frac{1}{2\mu}\!\|y\!-\!(Tx\!+\!\mu \lbb)\|^2\!-\!\frac{\mu}{2}\|\lbb\|^2.
\end{align*}
According to the proximal mapping defined in \eqref{p1}, the minimizer $y^{*}_{x, \lbb}$ of $L^{\mu}$ with respect to $y$ is
\begin{align}\label{ystar}
y^{*}_{x, \lbb}=\operatorname{prox}_{\mu \varphi}(Tx+\mu \lbb).
\end{align}
Evaluation of $L^{\mu}(x, y; \l, \lb, \lbb)$ at this minimizer gives the following PAL
\begin{align}\label{augl}
\mathscr L_{\mu}(x; \l, \lb, \lbb)=f(x)+\varphi_{\mu}(Tx+\mu \lbb)+\l^{\top} g(x)+\lb^{\top}h(x)-\frac{\mu}{2}\|\lbb\|^2.
\end{align}
By definition \eqref{mor}, this function is exactly the Moreau envelope of $L^{\mu}$   and thus is continuously differentiable with respect to both $x$ and $(\l, \lb, \lbb)$.

The PAL  was proposed in \cite{dhingra2019}, but without considering optimization constraints.
When such constraints are taken into account, a recent study \cite{wang2021} identifies three fundamental challenges in extending PAL to settings with non-smooth objectives, with exponential stability analysis emerging as a critical theoretical obstacle.
To prove exponential stability, the IQC framework \cite{lessard2016,fazlyab2018} provides systematic tools for exponential convergence analysis via linear matrix inequalities for systems  in feedback interconnections involving linear dynamics and IQC-satisfying nonlinearities.
While gradient and proximal mappings individually comply with IQC criteria, their compositions with nonlinear operators generally violate such constraints unless combined with affine operations \cite{fazlyab2018}. Consequently, existing IQC-based PAL analyses for exponential convergence predominantly address  equality constraints  \cite{dhingra2019,hassan2020,hassan2021,wang2021,ding2020} or affine inequalities \cite{hassan2018,bansode2019,qu2019}, leaving nonlinear convex inequality constraints unexplored.
While some existing works address general constraints, they either guarantee only asymptotic convergence \cite{goldsztajn2021} or achieve semi-global exponential convergence \cite{tang2020}. To handle general optimization constraints while ensuring exponential convergence, we bypass the IQC framework. Instead, we establish in the following lemma that the PAL $\mc L_{\mu}(x;\l, \lb, \lbb)$ is strongly convex in the primal variable $x$ and strongly concave in the dual variable $(\l, \lb, \lbb)$, a vital property for proving exponential convergence of our optimization algorithm.

\begin{lemma}\label{L4.1}
Under Assumptions \ref{sconvex} and \ref{lcontinuous}, the PAL defined in \eqref{augl} is $\a$-strongly convex in $x$ and $\left(\frac{\mu\ell}{\mu+\ell}+2\mu\right)$-strongly concave in $(\l, \lb, \lbb)$.
\end{lemma}

\noindent{\bf Proof:}   Obviously $L^{\mu}$  in \eqref{Lag1} is convex in $(x,y)$ since it is  a  linear combination of the convex and affine  functions  $f(x)$, $\varphi(y)$, $g(x)$, $h(x)$, $Tx-y$ and $\|Tx-y\|^2$ (the convexity of the  last one can be shown by a simple calculation) with the combination coefficients for convex functions being non-negative;  furthermore, $L^{\mu}$ is strongly convex in $x$ since $f(x)$ is strongly convex due to Assumption \ref{sconvex}.

We now show that  $\mathscr L_{\mu}(x; \l, \lb, \lbb)$ is strongly concave in $(\l, \lb, \lbb)$.
Note
\begin{align*}
\mathscr L_{\mu}(x; \l, \lb, \lbb)
&=f(x)+\l^{\top} g(x)+\lb^{\top}h(x)-\frac{\mu}{2}\|\lbb\|^2\\
&\hspace{2.5cm}-\sup_{y\in \R^m}\left\{-\varphi(y)-\frac{1}{2\mu}\|y-(Tx+\mu \lbb)\|^2\right\}\\
&=f(x)+\l^{\top} g(x)+\lb^{\top}h(x)-\frac{\mu}{2}\|\lbb\|^2-\frac{1}{2\mu}\|Tx+\mu\lbb\|^2\\
&\hspace{2.5cm}-\sup_{y\in \R^m}\left\{-\varphi(y)-\frac{1}{2\mu}\|y\|^2
+y^{\top}(\frac{1}{\mu}Tx+\lbb)\right\}\\
&\!=\!f(x)\!+\!\l^{\!\top}\! g(x)\!+\!\lb^{\!\top}\!h(x)\!-\!\frac{\mu}{2}\|\lbb\|^2\!-\!\frac{\mu}{2}\|\frac{Tx}{\mu}\!+\!\lbb\|^2\!-\!
\tilde\varphi^{*}(\frac{1}{\mu}Tx\!+\!\lbb),
\end{align*}
where $\tilde\varphi^*$ is the Fenchel conjugate of the function $\tilde\varphi(y)=\varphi(y)+\frac{1}{2\mu}\|y\|^2$.  Noting that $\tilde\varphi$ is $1/\mu$-strongly convex, it then follows from \cite[Theorem 4.2.2]{hiriart2004} that $\tilde\varphi^*$ is $\mu$-smooth.   Also, in view of the $1/\ell$-Lipschtz continuous of $\varphi$ assumed in Assumption \ref{lcontinuous} or equivalently the $(1/\ell+1/\mu)$-Lipschtz continuous of $\tilde \varphi$, one has  $\frac{\mu\ell}{\mu+\ell}$-strong convexity of $\varphi^*$. Therefore, $\tilde\varphi^{*}(\frac{1}{\mu}Tx+\lbb)$ is a $\frac{\mu\ell}{\mu+\ell}$-strongly convex function of $\lbb$. Also, both $\frac{\mu}{2}\|\lbb\|^2$ and  $\frac{\mu}{2}\|\frac{Tx}{\mu}+\lbb\|^2$ are $\mu$-strongly convex functions of $\lbb$. Therefore, $\mathscr L_{\mu}(x; \l, \lb, \lbb)$ is a $\left(\frac{\mu\ell}{\mu+\ell}+2\mu\right)$-strongly concave function of $(\l, \lb, \lbb)$.
\hfill $\blacksquare$

With the above lemma, we are directed toward utilizing $\mathscr L_{\mu}(x; \l, \lb, \lbb)$ as a Lagrangian to design corresponding primal-dual dynamics. A preliminary work should be done to relate  the optimal solution $x^*$ of \eqref{optimization} to the saddle point of $\mathscr L_{\mu}(x; \l, \lb, \lbb)$. This is shown in the following lemma.
\begin{lemma}\label{andian}
Consider the optimization problem \eqref{optimization} under Assumptions \ref{sconvex}-\ref{licq}.

{\rm (1)} Suppose there exists a point $x^*\in \R^n$ and multipliers $(\l^*, \lb^*, \lbb^*)\in \R^r_{+}\times \R^s \times \R^m$ such that
\begin{align}\label{saddle}
\mathscr L_{\mu}(x^*; \l, \lb, \lbb)\leq \mathscr L_{\mu}(x^*, \l^*; \lb^*, \lbb^*) \leq \mathscr L_{\mu}(x; \l^*, \lb^*, \lbb^*)
\end{align}
holds for all $(x; \l, \lb, \lbb)\in \R^n \times \R^r_{+}\times \R^s \times \R^m$. Then $x^*$ is an optimal solution to the optimization problem \eqref{optimization}.

{\rm (2)} Suppose  $x^*\in \R^n$ is an optimal solution to the optimization problem \eqref{optimization}. Then there exist  $(\l^*, \lb^*, \lbb^*)\in \R^r_{+}\times \R^s \times \R^m$ such that \eqref{saddle} holds.
\end{lemma}

{\bf Proof:}  (1) We show that $(x^*, y^*;\l^*, \lb^*, \lbb^*)$ with $y^*=y^{*}_{x, \lbb}$ defined in \eqref{ystar} satisfies the $\mc K\mc K\mc T$ system \eqref{KKT2} so that $(x^*, y^*)$ is an optimal solution to \eqref{optimization2} and consequently $x^*$ is an optimal solution to \eqref{optimization}.

Firstly, we use the first inequality in \eqref{saddle}, i.e.,
\begin{align*}
&f(x^*)+\varphi(y^{*}_{x^*, \lbb})+\l^{*\top} g(x^*)+\lb^{\top}h(x^*)+\lbb^{\top}(Tx^*-y^{*}_{x^*, \lbb})+\frac{1}{2\mu}\|Tx^*-y^{*}_{x, \lbb}\|^2\\
&\hfill \leq \mc L_{\mu}(x^*; \l^*, \lb^*, \lbb^*).
\end{align*}
Noting that the right side of the above inequality is a constant and this inequality holds for all $\l\succeq 0$, we must have $g(x^*)\preceq 0$. Similarity, this inequality holds for all $\lb$ so that one must have $h(x^*)=0$. Also, this inequality holds for all $\lb$, implying $Tx^*=y^{*}_{x, \lbb}$.

On the one hand, taking  $(\l, \lb)=(0, \lb^*)$ in the first inequality in \eqref{saddle} gives $\l^{*\top} g(x^*)\geq 0$; on the other hand, $\l^* \succeq 0$ and $g(x^*)\preceq 0$ give $\l^{*\top} g(x^*)\leq 0$. Therefore, $g(x^*)\preceq 0$ gives $\l^{*\top} g(x^*)=0$.

Lastly, we prove $(x^*;\l^*,\lb^*,\lbb^*)$ satisfies \eqref{KKT2f}.  With $y^*=y^{*}_{x^*, \lbb^*}$ and by referring to the right hand side of the inequality \eqref{saddle}, one has
\begin{align*}
 &L^{\mu}(x^*, y^{*}_{x^*, \lbb^*}; \l^*, \lb^*, \lbb^*)\\
 =&\mathscr L_{\mu}(x^*; \l^*, \lb^*, \lbb^*)\leq
 \mathscr L_{\mu}(x; \l^*, \lb^*, \lbb^*)\leq  L^{\mu}(x,  y^{*}_{x, \lbb^*}; \l^*, \lb^*, \lbb^*).
\end{align*}
This means that $x^*$ is a minimizer of $x\ra  L^{\mu}(x, y^{*}_{x, \lbb^*}; \l^*, \lb^*, \lbb^*)$ so that its subdifferential set includes $0$.  Noting that $ L^{\mu}(x^*, y; \l^*, \lb^*, \lbb^*)$ achieves it minimum at $y^{*}_{x^*, \lbb^*}$,  its derivative with respect to $y$ at $y^{*}_{x^*, \lbb^*}$ is zero. Now the necessary conditions for $x^*$ to be a minimizer of $x\ra  L^{\mu}(x, y^{*}_{x, \lbb^*}; \l^*, \lb^*, \lbb^*)$ becomes
\begin{align}\label{nec1}
\hspace{-0.3cm}0\!=\! \n f(x^*)\!+\!\sum_{i=1}^{r} \!\l^{*\!\top} \n g_i(x^*)\!+\!\sum_{i=1}^{s}\! \lb^{*\!\top}   \n h_i(x^*)\!+\!T^{\!\top}\lbb^{*}\!+\!\frac{1}{\mu}T^{\!\top} (Tx^*\!-\!y^{*}_{x^*\!\!, \lbb^*})
\end{align}
Also noting that $y^{*}_{x^*, \lbb^*}$ is a minimizer of the function $y\ra  L^{\mu}(x^*, y;\l^*, \lb^*, \lbb^*)$, the following optimality condition holds
\begin{align}\label{nec2}
\partial \varphi(y^{*}_{x^*, \lbb^*})-\lbb^{*}-\frac{1}{\mu} (Tx^*-y^{*}_{x^*, \lbb^*}) \ni 0
\end{align}
These two results \eqref{nec1} and \eqref{nec2}, together with $y^{*}_{x^*, \lbb^*}=Tx^*$,  give
\begin{align*}
\n f(x^*)+\sum_{i=1}^{r} \l^{*\top} \n g_i(x^*)+\sum_{i=1}^{s} \lb^{*\top}   \n h_i(x^*)+T^{\top}\lbb^*p\ni 0.
\end{align*}
Therefore, $(x^*; \l^*, \lb^*, \lbb^*)$ satisfies the KKT system \eqref{KKT2}, implying that  $(x^*, y^*)$ is an optimal solution  of \eqref{optimization2} and consequently $x^*$ is an optimal solution of \eqref{optimization}.

(2) Suppose  $x^*\in \R^n$ is an optimal solution to the optimization problem \eqref{optimization} so that $(x^*, y^*)$ with $y^*=Tx^*$ is an optimal solution to \eqref{optimization2}, implying the existence of $(\l^*,\lb^*,\lbb^*)$ such that $(x^*, y^*;\l^*,\lb^*,\lbb^*)$ satisfies  the KKT system \eqref{KKT2}.  Therefore,
\begin{align}\label{eq1}
\mathscr L_{\mu}(x^*;\l^*,\lb^*,\lbb^*)\leq  L^{\mu} (x^*, y^*;\l^*,\lb^*,\lbb^*)=f(x^*)+\varphi(y^*).
\end{align}
On the other hand, the strong duality for the problem \eqref{optimization2} implies
\begin{align}\label{eq2}
f(x^*)+\varphi(y^*)
&\!=\!\underset{x,y}\min\left\{ f(x)+\varphi(y)+\l^{*\top}g(x)+\lb^{*\top}h(x)+\lbb^{*\top}(Tx-y) \right\}\nonumber\\
&\leq f(x)+\varphi(y)+\l^{*\top}g(x)+\lb^{*\top}h(x)+\lbb^{*\top}(Tx-y) \nonumber \\
& \leq \mathscr L_{\mu} (x; \l^*,\lb^*,\lbb^*).
\end{align}
Combing the inequalities \eqref{eq1} and \eqref{eq2}, one obtains  $\mathscr L_{\mu}(x^*;\l^*,\lb^*,\lbb^*) \leq \mathscr L_{\mu} (x; \l^*,\lb^*,\lbb^*)$ which is the second inequality in \eqref{saddle}.

To prove the first inequality in \eqref{saddle},  we recall the KKT systems \eqref{KKT2} and have the following inequality
\begin{align*}
\mathscr L_{\mu} (x^*; \l,\lb,\lbb)
&=f(x^*)+\varphi(y^{*}_{x^*, \lbb^*})+\l^{\top} g(x^*)+\lb^{\top}h(x^*)+\lbb^{\top}(Tx^*-y^{*}_{x^*, \lbb^*})\\
&\hspace{5.7cm}+(1/2\mu)\|Tx^*-y^{*}_{x^*, \lbb^*}\|^2\\
&=f(x^*)+\varphi(y^{*}_{x^*, \lbb^*})\!+\!\l^{*\top} g(x^*)\!+\!\lb^{*\top}h(x^*)\!+\!\lbb^{*\top}(Tx^*\!-\!y^{*}_{x^*\!\!, \lbb^*})\\
&\hspace{5.7cm}+(1/2\mu)\|Tx^*-y^{*}_{x^*, \lbb^*}\|^2\\
& =\mathscr L_{\mu}(x^*; \l^*,\lb^*,\lbb^*).
\end{align*}
This completes the second part of the proof.
\hfill $\blacksquare$

\subsection{Projection-free primal-dual optimization dynamics}

We are now in a position to design a primal-dual dynamics building on the PAL \eqref{augl} to seek the optimal solution $(x^*, y^*)$ of the optimization problem \eqref{optimization2} with exponential convergence.  This algorithm excuses gradient descent on $x$ and gradient ascent on $(\l, \lb, \lbb)$.   However, the gradient descent of $\mathscr L_{\mu}$ on $\l$ requires a projection onto the positive quadrant to keep the resulting $\l$ dynamics stay non-negative; more specifically, $\dot \l= [\n_{\l}\mathscr L_{\mu}]_{+}$, which has a discontinuous right-hand side. In this paper, we propose a new design procedure for this Lagrangian dynamics without using projection but referring to mirror ascent; it comes with the advantage that the resulting dynamics is smooth so that its convergence analysis and numerical simulations  are relatively simple compared with non-smooth projection based dynamics. More specifically, we propose the following projection-free algorithm which, for the PAL \eqref{augl}, performs gradient descent on the primal variable $x$, gradient ascent on the dual variables $(\lb, \lbb)$, and mirror ascent on the dual variable $\l$, specified as follows,
\begin{subequations} \label{proxPD}
    \begin{empheq}[left={ \empheqlbrace\,}]{align}
&\dot x=-\n f(x)-T^{\top}\n \varphi_{\mu}(Tx+\mu \lbb)-\l^{\top} \n g(x)-\lb^{\top}\n h(x), \label{sf1}\\
&\dot \l=[\l \oslash (1+\eta \odot \l)] \odot g(x),  \l(0) \succeq 0, \label{sf2}\\
&\dot \lb=h(x), \label{sf3}\\
&\dot \lbb=\mu \n \varphi_{\mu}(Tx+\mu \lbb)-\mu \lbb. \label{sf4}
    \end{empheq}
\end{subequations}

\begin{remark}
Some remarks are given regarding the intuition behind the projection-free dynamics \eqref{sf2}.
Recall that the mirror descent algorithm is devoted to the constrained minimization problem $\min_{y \in \mc Y}\digamma(y)$, where
 $\mc Y$ is a convex  set in an Euclidian space. Let $\phi: \mc Y\!\ra\! \mb{R}$ be a function which is strictly convex and twice differentiable.  Its conjugate convex function
$\phi^*: \mc Y\!\ra\! \mb{R}$  is defined as
$
\phi^*(\omega)\!=\!\sup_{y\in \mc Y}\{(y, \omega)\!-\!\phi(y)\},
$
which can be shown as  convex and  differentiable \cite[Theorem 26.3]{rockafellar1970}.
Let $\mc Z\!=\!\{z|z\!=\!\nabla \phi(y), y \!\in\! \mc Y\}$ be the image of $\mc Y$ under the mapping $\nabla \phi$. Then $\nabla \phi: \mc Y\!\ra\! \mc Z$ and  $\nabla \phi^*: \mc Z\!\ra\! \mc Y$. The continuous-time mirror descent algorithm for the above constrained minimization problem takes the following form (\cite[e.q. (5)]{raginsky2012}): $\dot z\!=\!-\!\n \digamma(y), y\!=\!\n \phi^*(z)$.
The second equation $ y\!=\!\n \phi^*(z)$ is equivalent to $z\!=\!\n \phi(y)$. Taking derivative yields $\dot z\!=\!\n^2 \phi(y) \dot y$. Therefore the first equation in the mirror descent becomes $\dot y\!=\!-[\n^2 \phi(y)]^{-1}\n \digamma(y)$.
 Choosing the constraint set as $\mc Y\!=\!\{y|y_i \geq 0\}$ and $\phi(y)\!=\!\frac{\eta}{2} \|y\|^2\!+\!\sum_{i=1}^{n}y_i\ln y_i$ which is well defined on $\mc Y$,   the above mirror dynamics can be written as
$\dot y=-{\rm diag}[y_i/(1+\eta y_i)]_{i=1}^n \nabla \digamma(y)$.
Letting $y=\l_{i}$ and $\digamma(\l_{i})$ be a function of  $\l_{i}$ defined as the right hand side of equation \eqref{augl}, then the above dynamics is exactly \eqref{sf2}.
\end{remark}

Before carrying out the convergence analysis of the algorithm \eqref{proxPD}, the relationship between the equilibrium of \eqref{proxPD} and the optimal solution of \eqref{optimization} is revealed in the following lemma.
\begin{lemma}\label{equilibrium}
A point $(x^*; \l^*, \lb^*, \lbb^*)\in \R^n\times\R^r_{+}\times \R^s \times \R^m$  is the equilibrium of \eqref{proxPD} if and only if  $x^*$ is the optimal solution to the optimization problem \eqref{optimization2}.
\end{lemma}

\noindent{\bf Proof:}  Suppose $x^*$ is an optimal solution. It follows from Lemma \ref{andian} that there exists $(\l^*, \lb^*, \lbb^*)\in \R^r_{+}\times \R^s \times \R^m$ such that the saddle point inequality \eqref{saddle} holds. This implies that at $(\l^*, \lb^*, \lbb^*)$ it holds $\partial \mathscr{L}_{\mu}(x^*; \l^*, \lb^*, \lbb^*)/\partial x=0$, $\partial \mathscr{L}_{\mu}(x^*; \l^*, \lb^*, \lbb^*)/\partial \lb=0$ and $\partial \mathscr{L}_{\mu}(x^*; \l^*, \lb^*, \lbb^*)/\partial \lbb=0$. These three equations exactly shows that the right hands of the equations \eqref{sf1}, \eqref{sf3}, and \eqref{sf4} vanish at $(\l^*, \lb^*, \lbb^*)$.  The right hand of the equation \eqref{sf2} also takes zero value at $(\l^*, \lb^*, \lbb^*)$ in view of KKT subsystem \eqref{KKT2e}. Therefore, $(x^*; \l^*, \lb^*, \lbb^*)$  is the equilibrium of \eqref{proxPD}.

On the other hand, letting $(x^*; \l^*, \lb^*, \lbb^*)$ be the equilibrium of \eqref{proxPD} and inserting it into \eqref{sf1}, \eqref{sf3}, and \eqref{sf4}, one obtains
\begin{align}\label{partial1}
&\frac{\partial \mathscr L_{\mu}(x^*; \l^*, \lb^*, \lbb^*)}{\partial x}=0,\\
&\frac{\partial \mathscr L_{\mu}(x^*; \l^*, \lb^*, \lbb^*)}{\partial \lb}=0, \label{partial2}\\
&\frac{\partial \mathscr L_{\mu}(x^*; \l^*, \lb^*, \lbb^*)}{\partial \lbb}=0.  \label{partial3}
\end{align}
Recall the fact that, for a differentiable convex (concave) function $F(x)$, a necessary and sufficient condition for $\bar x$ to be the minima (maxima) of $F$ is $\n F(\bar x)=0$.
By using this fact,  the result  \eqref{partial1} gives the minima inequality
\begin{align}\label{zuobian}
\mathscr L_{\mu}(x^*; \l^*, \lb^*, \lbb^*)\leq \mathscr L_{\mu}(x, \l^*, \lb^*, \lbb^*),
\end{align}
and the results in  \eqref{partial2}-\eqref{partial3} indicates  the maxima inequality
$\mathscr L_{\mu}(x^*;  \l^*, \lb, \lbb)$ $\leq$ $\mathscr L_{\mu}(x^*;  \l^*, \lb^*, \lbb^*)$
with respect to $(\lb, \lbb)$. Due to the decoupling of $\l$ with $(\lb, \lbb)$ in the PAL \eqref{augl}, this inequality holds for all $\l$ rather than only on $\l^*$,
\begin{align}\label{t1}
\mathscr L_{\mu}(x^*;   \l, \lb, \lbb) \leq \mathscr L_{\mu}(x^*;  \l, \lb^*, \lbb^*).
\end{align}
Inserting $(x^*;  \l^*, \lb^*, \lbb^*)$ into  \eqref{sf2} give $\l^* \odot g(x^*)=0$; i.e., $\l_i^*g_i(x^*)=0, i=1, \cdots, r$. This shows that $g_i(x^*)=0$ or $\l_i^*=0$. The former is nothing but $\partial \mathscr L_{\mu}(x^*;  \l^*, \lb^*, \lbb^*)/\partial \l=0$ so that the function $\l\ra \mathscr L_{\mu}(x^*;  \l, \lb^*, \lbb^*)$ is maximized  at $\l^*$; the latter, together with $g_i(x^*)\leq 0$,  shows that the function $\l\ra \mathscr L_{\mu}(x^*;  \l, \lb^*, \lbb^*)$ is maximized  at $\l^*$. Both cases indicate that
\begin{align}\label{t2}
\mathscr L_{\mu}(x^*;  \l, \lb^*, \lbb^*) \leq \mathscr L_{\mu}(x^*; \l^*,  \lb^*, \lbb^*).
\end{align}
The combination of \eqref{t1} and \eqref{t2} gives
\begin{align}\label{youbian}
\mathscr L_{\mu}(x^*;  \l,  \lb,  \lbb) \leq \mathscr L_{\mu}(x^*; \l^*,  \lb^*, \lbb^*).
\end{align}
The two  inequalities \eqref{zuobian} and \eqref{youbian} tell us that $(x^*; \l^*,  \lb^*, \lbb^*)$ is the saddle point of $\mathscr L_{\mu}$.  According to the first result in Lemma \ref{andian}, $x^*$ is an optimal solution to the optimization problem \eqref{optimization2}.  \hfill $\blacksquare$

With the relationship between the equilibrium of the optimization dynamics \eqref{proxPD} and the optimal solution to the problem \eqref{optimization2} discovered in the above lemma, we have the following result which proves the exponential stability of  our algorithm \eqref{proxPD}.
\begin{theorem}\label{expSta}
For the optimization problem \eqref{optimization2} under Assumptions \ref{sconvex}-\ref{licq}, consider the optimization algorithm \eqref{proxPD}. Then, for any trajectory of \eqref{proxPD} with initial conditions $x(0)\in \R^n$ and $(\l(0), \lb(0), \lbb(0))\in \R^r_{+}\times \R^s \times \R^m$, the sub-trajectory $x(t)$ converges to the optimal solution $x^*$ exponentially as $t\ra +\w$.
\end{theorem}

\noindent{\bf Proof:} Construct a Lyapunov candidate  $V(x; \l, \lb, \lbb)=V_1+V_2+V_3+V_4$ with
\begin{align*}
&  V_1=\frac12 \|x-x^*\|^2,\\
&  V_2=\frac12 \sum\nolimits_{i=1}^{r}\eta_{i}(\lambda_{i}-\lambda_{i}^*)^2,\\
&  V_3=\sum\nolimits_{i\in \Omega}D_{\psi}(\lambda_{i}, \lambda_{i}^*) +\sum\nolimits_{i\notin \Omega}(\lambda_{i}-\lambda_{i}^*)^2,\\
 & V_4=\frac12  \sum\nolimits_{i=1}^{s}(\lb_{i}-\lb_{i}^*)^2+\frac12  \sum\nolimits_{i=1}^{s}(\lbb_{i}-\lbb_{i}^*)^2,
\end{align*}
where $\Omega=\{i| \lambda_{i}^* \!\neq\! 0\}$ and $D_{\psi}(\lambda_{i}, \lambda_{i}^*)\geq 0$ is the Bregman divergence \cite{bregman1967} between $\lambda_{i}$ and $\lambda_{i}^*$ with respect to $\psi(t)=t\ln t$.
Therefore, $V\geq 0$,  and also $V(x; \l, \lb, \lbb)=0$ if and only if $(x; \l, \lb, \lbb)\!=\!(x^*; \l^*, \lb^*, \lbb^*)$.  We now use this candidate Lyapunov function to prove exponential stability of the equilibrium $(x^*; \l^*, \lb^*, \lbb^*)$ of the system \eqref{proxPD} by referring to the criterion in \cite[Corollary 3.4]{khalil1996}. To this end, we will show that $V$ is lower and upper bounded by quadratic functions and the time derivative of $V$ is upper bounded by a negative quadratic function. The following two parts are the details.

(1) We first show that $V$ can be upper and lower bounded by quadratic functions of $(x; \l, \lb, \lbb)$ (we only prove this for $V_3$ since this is obvious for $V_1, V_2, V_4$).
For $i\in \O$ which implies $\l_i^*>0$,  the Bregman $D_{\psi}(\l_i, \l_i^*)$ can be calculated as
\begin{align*}
D_{\psi}(\l_i, \l_i^*)=\psi(\l_i)-\psi(\l_i^*)-\psi'(\l_i^*)(\l_i-\l_i^*)=\frac{1}{2\varsigma_i}(\l_i-\l_i^*)^2.
\end{align*}
with $\varsigma_i$ lying between $\l_i^*$ and $\l_i$. In view of $\l_i^*>0$ and  $\l_i\ra \l_i^*$, $\varsigma_i$ is a positive constant. Denoting $\tilde a=\underset{i\in \O}{\min}\{1/2\varsigma_i\}$ and  $\tilde b=\underset{i\in \O}{\max}\{1/2\varsigma_i\}$, one has $\tilde a\leq D_{\psi}(\l_i, \l_i^*)\leq \tilde b$. With this, one can obtain the following lower and upper quadratic bounds for $V_3$,
\begin{align*}
a\|\l-\l^*\|^2
&\leq \sum\limits_{j\in \Omega}\tilde a |\l_i-\l_i^*|^2 +\sum\limits_{j\notin \Omega}(\lambda_{j}-\lambda_{j}^*)^2  \\
&\leq V_3  \leq \sum\limits_{j\in \Omega}\tilde b |\l_i-\l_i^*|^2 +\sum\limits_{j\notin \Omega}(\lambda_{j}-\lambda_{j}^*)^2 \leq b\|\l-\l^*\|^2,
\end{align*}
where  $a=\min(\tilde a, 1)>0$ and  $b=\max(\tilde b, 1)>0$.   Therefore, the function $V$ has
quadratic lower and upper bounds.

(2) The time derivative of $V$ along the trajectories of \eqref{proxPD} is,
\begin{align*}
\dot V|_{(\ref{proxPD})}
&=-(x-x^*)^{\top} \frac{\partial  \mathscr L_{\mu}}{\partial x}+\sum_{i=1}^{r} \frac{\eta_{i}\lambda_{i}(\lambda_{i}-\lambda_{i}^*)}{1+\eta_{i}\lambda_{i}}  \frac{\partial  \mathscr L_{\mu}}{\partial \l}+\sum_{i=1}^{r}  \frac{\lambda_i}{1+\eta_i\lambda_i}  \frac{\partial \mathscr L_{\mu}}{\partial \l}\\
&\hspace{0.5cm}-\sum_{i\in \Omega} \frac{\lambda_i^*}{1+\eta_i\lambda_i}  \frac{\partial  \mathscr L_{\mu}}{\partial \l}
+\sum_{i=1}^{s} (\lb_i-\lb_i^*) \frac{\partial \bar L}{\partial \lb}+\sum_{i=1}^{m} (\lbb_i-\lbb_i^*) \frac{\partial  L_{\mu}}{\partial \lbb}\\
&=-(x\!-\!x^*)^{\top} \frac{\partial   \mathscr L_{\mu}}{\partial x}\!+\! (\lambda\!-\!\lambda^*)^{\top} \frac{\partial  \mathscr L_{\mu}}{\partial \l}\!+\! (\lb\!-\!\lb^*)^{\top} \frac{\partial   \mathscr L_{\mu}}{\partial \lb}\!+\! (\lbb\!-\!\lbb^*)^{\top} \frac{\partial   \mathscr L_{\mu}}{\partial \lbb}\\
&=-(x-x^*)^{\top} \frac{\partial   \mathscr L_{\mu}}{\partial x}+ (\L-\L^*)^{\top} \frac{\partial \mathscr L_{\mu}}{\partial \L},
\end{align*}
where $\L=(\l, \lb, \lbb)$ and $\L^*=(\l^*, \lb^*, \lbb^*)$.
Noting that $\mathscr L_{\mu}(x; \l, \lb, \lbb)$ has been proved in Lemma \ref{L4.1} to be  strongly convex in $x$ and strongly concave in $(\l, \lb, \lbb)$,  the above inequality becomes
\begin{align}\label{dotV}
\dot V|_{(\ref{proxPD})}
&\leq [ \mathscr L_{\mu}(x^*; \l, \lb,\lbb)- \mathscr L_{\mu}(x; \l, \lb,\lbb)]+[ \mathscr L_{\mu}(x; \l, \lb,\lbb)- \mathscr L_{\mu}(x; \l^*, \lb^*,\lbb^*)]\nonumber\\
&\hspace{0.5cm} -\a \|x-x^*\|^2-\left(\frac{\mu\ell}{\mu+\ell}+2\mu\right)\|\L-\L^*\|^2.
\end{align}
This shows that $\dot V$ has minus quadratic upper bound.
\hfill $\blacksquare$

Usually, in designing algorithms for constrained optimization
problem, either hard constraints or soft constraints can be considered. In the former, the
optimization constraints should be satisfied on the fly (e.g., the interior algorithms),
while in the latter, rather than forcing each transient state of the iterate to obey the constraints, one requires the optimization constraints to be satisfied only at the asymptotic
state. Since the ultimate objective is to compute the optimal solution that satisfies the constraints, our algorithm just follows the second line of soft constraints and therefore there is
no need in this paper to show positive feasibility for each transient state $\l(t)$ but feasibility for the asymptotic state of $\l(t)$ only. This is indeed the case in view of the equilibrium analysis in Lemma \ref{equilibrium}.

\begin{remark} We note that our method is capable of achieving fully smooth dynamics while handling general convex inequality constraints.
In primal-dual optimization dynamics based on the PAL framework,  the multiplier $\l$ associated with inequality constraint must remain nonnegative. Existing methods enforce this by incorporating a projection onto the positive quadrant \cite{feijer2010,goldsztajn2021,tang2020}, but this introduces nonsmooth dynamics that complicate stability analysis.
To address this, smooth alternatives have been proposed to eliminate projection.
The works \cite{qu2019,adegbege2021,bin2024} handle inequality constraints without explicit projection, instead they replace projection-based updates with max-based regularization (specifically algorithm (9a)-(9c) in \cite{adegbege2021}, (9a)-(9b) in \cite{qu2019}, and (1b) in \cite{bin2024}), implicitly ensuring nonnegativity of the multiplier associated with inequality constraint. While this improves computational simplicity, the max operator still makes the optimization dynamics nonsmooth in the classical sense.
In contrast, our  projection-free dynamics \eqref{sf2} is fully smooth, and this smooth nature enables rigorous stability certification via classical Lyapunov techniques, bypassing the need for nonsmooth analysis tools like differential inclusions.
More important, the max-based method in \cite{qu2019,adegbege2021}  are limited to affine inequality constraints, with no extension to general convex inequality constraints. Our approach overcomes this limitation and  extends to general convex inequality constraints.
\end{remark}

\begin{remark}
Asides from its ability to obtain fully smooth optimization dynamics and deal with general convex inequality constraints, our method can yield exponential convergence rather than asymptotic convergence.
As  noted in \cite{wang2021}, the combination of general convex constraints and nonsmooth objectives poses a significant challenge to achieving exponential stability. For instance, existing studies on general inequality constraints typically ensure only asymptotic convergence \cite{goldsztajn2021} or semi-global exponential convergence \cite{tang2020}. IQC is a widely used framework for exponential convergence analysis.
However, existing IQC-based approaches primarily address affine-constrained settings, such as equality constraints \cite{dhingra2019,hassan2020,hassan2021,wang2021,ding2020} or affine linear inequalities \cite{hassan2018,bansode2019}, leaving nonlinear convex inequality constraints unexamined. Similar limitations appear in projection-free methods \cite{adegbege2021} and exponential stability analyses \cite{qu2019}, which remain restricted to affine inequalities.
We do not use the IQC method, we can prove exponential stability without assuming affine linearity assumption on inequality constraints.
\end{remark}

Note that, in the proof of Theorem \ref{expSta}, the strong concavity of the PAL function $\mathscr L_{\mu}$ on $(\l, \lb, \lbb)$ depends only on  the Lipschitz continuous of $\varphi$  in the sense of Assumption \ref{lcontinuous}, but not on the strong convexity of $f$.
If  the strong convexity assumption on $f$ made in Assumption \ref{sconvex} is weaken to be convex, the PAL $\mathscr L_{\mu}$ in \eqref{augl} is only  convex on $x$, but still strongly concave on $(\l, \lb, \lbb)$. In this case, only asymptotic stability of the optimization dynamics \eqref{proxPD} can be ensured. Before presenting the result in this case, the following assumption is given.

\begin{assbis}\label{convex}
The function $f$ is  convex and continuously differentiable.
\end{assbis}

\begin{theorem}
For the optimization problem \eqref{optimization2} under Assumptions \ref{convex}, \ref{lcontinuous} and \ref{licq}, consider the optimization algorithm \eqref{proxPD}. Then, for any trajectory of \eqref{proxPD} with initial conditions $x(0)\in\R^n$ and $(\l(0), \lb(0), \lbb(0))\in \R^r_{+}\times \R^s \times \R^m$, the sub-trajectory $x(t)$ converges to the optimal solution $x^*$ asymptotically as $t\ra +\w$.
\end{theorem}

\noindent{\bf Proof:}  Consider also the Lyapunov function constructed in the proof  of Theorem \ref{expSta}. The time derivative of $V$ calculated in \eqref{dotV} is now modified as
\begin{align}\label{dotV'}
\dot V|_{(\ref{proxPD})}
&\leq [ \mathscr L_{\mu}(x^*; \l, \lb,\lbb)- \mathscr L_{\mu}(x; \l, \lb,\lbb)] \nonumber\\
&+[ \mathscr L_{\mu}(x; \l, \lb,\lbb)- \mathscr L_{\mu}(x; \l^*, \lb^*,\lbb^*)]-\left(\frac{\mu\ell}{\mu+\ell}+2\mu\right)\|\L-\L^*\|^2. \nonumber\\
&\leq [ \mathscr L_{\mu}(x^*; \l, \lb,\lbb)- \mathscr L_{\mu}(x^*; \l^*, \lb^*,\lbb^*)]\nonumber\\
&+[ \mathscr L_{\mu}(x^*; \l^*, \lb^*,\lbb^*)- \mathscr L_{\mu}(x; \l^*, \lb^*,\lbb^*)]-\left(\frac{\mu\ell}{\mu+\ell}+2\mu\right)\|\L-\L^*\|^2.
\end{align}
Note that the $x$-related quadratic form in the right hand side of \eqref{dotV} is dropped due to the lack of strong convexity of $\mathscr L_{\mu}$. Since $(x^*; \l^*, \lb^*,\lbb^*)$  is a saddle point of $\mathscr L_{\mu}$, the two curl brackets are both nonnegative. Therefore, the right hand side of \eqref{dotV'} is less than or equal to zero. Thus $\dot V\leq 0$.

We now apply the LaSalle's invariance principle to prove the asymptotic stability; that is, we only need to prove that $\dot V=0$ implies $(x;\l,\lb,\lbb)=(x^*;\l^*, \lb^*, \lbb^*)$.  To this end, letting $\dot V=0$ gives   $\L=\L^*$  and $\mathscr L_{\mu}(x; \l^*, \lb^*,\lbb^*)=\mathscr L_{\mu}(x^*; \l^*, \lb^*,\lbb^*)$. Noting again $\L=\L^*$, one has
\begin{align*}
\mathscr L_{\mu}(x; \l, \lb,\lbb)
=\mathscr L_{\mu}(x; \l^*, \lb^*,\lbb^*)=\mathscr L_{\mu}(x^*; \l^*, \lb^*,\lbb^*)\leq \mathscr L_{\mu}(x^*; \l^*, \lb^*,\lbb^*),
\end{align*}
or $\mathscr L_{\mu}(x; \l, \lb,\lbb)\leq \mathscr L_{\mu}(x; \l^*, \lb^*,\lbb^*)\leq \mathscr L_{\mu}(x^*; \l^*, \lb^*,\lbb^*)$.
This shows that $(x; \l^*, \lb^*,\lbb^*)$ is a saddle point of $\mathscr L_{\mu}$ and therefore $x$ is the optimal solution of the problem \eqref{optimization2} in view of the first result in Lemma \ref{andian}. By uniqueness of the optimal solution, we have $x=x^*$. Thus, $(x;\l, \lb, \lbb)=(x^*; \l^*, \lb^*, \lbb^*)$. As a consequence of LaSalle's invariance principle, the equilibrium  $(x^*; \l^*, \lb^*, \lbb^*)$ is asymptotically stable. \hfill $\blacksquare$

\begin{remark}
Although small values of the regularization parameter $\mu$ would cause slow convergence, this slowness can be alleviated via the ``continuation schemes'' carried out in the following way. Firstly, given a particular parameter $\mu$, one runs the algorithm to obtain a solution. With this solution as the ``warm-start'' and taking it as the initial state of the second round of execution of the optimization dynamics with a smaller parameter $\mu$, one obtains a second solution.  Continuing this procedure with a sequence of decreasing $\mu$, the optimal solution can be computed in succession. This kind of ``continuation schemes'' have been  found quite effective in speeding up the algorithm \cite{figueiredo2007,afonso2010}.
\end{remark}


\section{Simulation example}\label{example}

As an example, we applied our algorithm to the following distributed optimization problem on a network with $N$ nodes,
\begin{eqnarray}\label{optimization3}
\left\{
\begin{array}{llll}
  {\rm minimize}   && \sum_{i=1}^{N}f_i(x), \label{objective}\\
  {\rm subject~to} &&  g_i(x)\preccurlyeq 0,  \label{constraint}\\
  {\rm }           && h_i(x)=0, i=1, \cdots, N. \label{constraint2}
\end{array}
\right.
\end{eqnarray}
Here $x\in \R^n$ is the decision variable,   $f_i: \mathbb{R}^n \rightarrow \mathbb{R}$ is the local cost function on node $i$,
$g_i=(g_{i1}, \cdots, g_{ir_i})^{\top}: \mathbb{R}^n \rightarrow \mathbb{R}^{r_i}$ is the inequality constrained function on node $i$, and
$h_i=(h_{i1}, \cdots, h_{is_i})^{\top}:\mathbb{R}^n \rightarrow \mathbb{R}^{s_i}$   is the equality constrained function on node $i$. If there is no constraints  for agent $i$, one simply sets corresponding constraint functions to be zero.

In view of large $N$, it is difficult for a single agent to do the optimization task. Instead, one usually uses $N$ agents to cooperatively solve the problem, where each agent $i$ executes a subtask of local optimization by only using information from agent $i$ (i.e., the information $f_i, g_i, h_i$ as well as their gradients) and those information from its neighbors $\mc N_i$.  The neighboring relationship encodes the cooperation among these agents and  is described by a graph. The cooperation makes local optimization algorithms compute the optimal solution in a consensus way. That is, letting $x_i\in\R^n$ denote the estimation of the optimal solution $x^*$ by local algorithm of agent $i \in \{1, \cdots, N\}$, the cooperation among these agents should be directed toward rendering $x_1=x_2=\cdots=x_N=x^*$.
Denoting respectively the graph Laplacian and incidence matrix of the graph by $L$ and $T$ so that $L=T^{\top}T$, the  consensus  $x_1=x_2=\cdots=x_N$ among agents can be enforce by imposing the equality constraint $(T\otimes I_n)\hx=0$ with $\hx={\rm col}(x_1, \cdots, x_n)\in \R^{nN}$ if the graph is connected.
By adding to the objective function in \eqref{optimization3} the following penalty function
\begin{align}\label{varphi2}
\varphi(T\hx)=
\left\{
\begin{array}{lllll}
&  0, &T\hx=0 \\
&  +\w, &\text{otherwise},
\end{array}
\right.
\end{align}
the optimization problem \eqref{optimization3} can be represented as the following form
\begin{eqnarray}\label{optimization4}
\left\{
\begin{array}{llll}
  {\rm minimize}   && \sum_{i=1}^{N}f_i(x_i)+\varphi(T\hx), \label{objective}\\
  {\rm subject~to} &&  g_i(x_i)\preccurlyeq 0,  \label{constraint}\\
  {\rm }           && h_i(x_i)=0, i=1, \cdots, N. \label{constraint2}
\end{array}
\right.
\end{eqnarray}
Defining the objective function $F(\hx)=\sum_{i=1}^{N}f_i(x_i)$ and stack constraint functions $G(\hx)={\rm col}[g_1(x_1), \cdots, g_N(x_N)]$, $H(\hx)={\rm col}[h_1(x_1), \cdots, h_N(x_N)]$, the problem \eqref{optimization4} can be written in a compact form as
\begin{eqnarray}\label{optimization5}
\left\{
\begin{array}{llll}
  {\rm minimize}   && F(\hx)+\varphi(T\hx), \label{objective}\\
  {\rm subject~to} &&  G(\hx)\preccurlyeq 0,  \label{constraint}\\
  {\rm }           && H(\hx)=0. \label{constraint2}
\end{array}
\right.
\end{eqnarray}

We now apply the algorithm \eqref{proxPD} to the nonsmooth composite  optimization problem \eqref{optimization5}. Note that, for the function $\varphi$ defined in \eqref{varphi2}, it is easy to show that $\n \varphi(y)=\frac{1}{\mu} y$.
Furthermore, corresponding to the multipliers $\l$ and $\lb$ in \eqref{proxPD} for a single agent, we use $\l_i=(\l_{i1}, \cdots, \l_{ir_i})^{\top}$ and $\lb_i=(\lb_{i1}, \cdots, \lb_{is_i})^{\top}$ to respectively denote multipliers corresponding to constraints $g_i(x)\preceq 0$ and $h_i(x)=0$; more specifically $\l_{ij}$ is the multiplier for $g_{ij}(x)\leq 0$ for $j=1,\cdots, r_i$ and $\lb_{ij}$ is the multiplier for $h_{ij}(x)= 0$ for $j=1,\cdots, s_i$. Note that if there is no inequality or equality constraint for agent $i$, we simply set $\l_i$ or $\lb_i$ to be zero vector of appropriate dimension.
As for the multiplier $\hlbb$ (bold lambda used here to account for multiple agents) for the constraints $T\hx=\hy$ in \eqref{proxPD}, we introduce a new multiplier $\hlbb'=T\hlbb$ with partition $\hlbb'={\rm col}[\lbb_1', \cdots, \lbb_N']$. With this new multiplier, the subequatons \eqref{sf1} and \eqref{sf4} have a simpler forms. More specifically, the resulting optimization algorithm has the following distributed form
\begin{subequations} \label{algorithm3}
 \begin{empheq}[left={\empheqlbrace\,}]{align}
  &\dot x_i\!=\!\frac{1}{\mu}\!\sum_{j\in \mc{N}_i}\!(x_j\!-\!x_i)\!-\!\nabla f_i(x_i) \!-\!\sum_{j}^{r_i}\!\lambda_{ij}\nabla g_{ij}(x_i)\!-\!\sum_{j}^{s_i}\!\lb_{ij}\nabla h_{ij}(x_i)\!-\!\lbb_i'\\
 &\dot \l_{ij}=\frac{\l_{ij}}{1+\eta_{ij}\l_{ij}}  g_{ij}(x_i), \hspace{0.5cm}j=1, \cdots, r_i, \\
&\dot \lb_{ij}=h_{ij}(x_i),\hspace{2cm} j=1, \cdots, s_i, \\
&\dot \lbb_i'=-\sum_{j\in \mc{N}_i}(x_j-x_i),
    \end{empheq}
\end{subequations}
where $\mc N_i=\{j|(i,j)\in \mathcal E\}$ is the neighboring agents of agent $i$.  This algorithm is distributed since the optimization dynamics for agent $i$ depends on information from its neighboring agents. See also \cite{dhingra2019} for  distributed implementation without considering optimization constraints.

\begin{figure}[htbp]
\centering{\includegraphics[width=2cm]{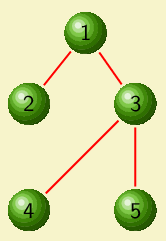}}
\caption{The topology of five agents described by a graph.}
\label{graph}
\end{figure}

We now apply the above distributed optimization algorithm  to  the following Rosen-Suzuki problem \cite{charalamous1977}
\begin{eqnarray*}
 & \min & x_1^2+x_2^2+2x_3^2+x_4^2-5x_1-5x_2-21x_3+7x_4\\
 & s.t. & -8+x_1-x_2+x_3-x_4+x_1^2+x_2^2+x_3^2+x_4^2\leq 0,\\
 &  & -10-x_1-x_4+x_1^2+2x_2^2+x_3^2+2x_4^2\leq 0, \\
 &  & -5+2x_1-x_2-x_4+2x_1^2+x_2^2+x_3^2=0.
\end{eqnarray*}
The optimal solution is shown in \cite{charalamous1977} to be $x^*=(0, 1, 2, -1)$.

\begin{figure}[htbp]
\centering{\includegraphics[width=12 cm]{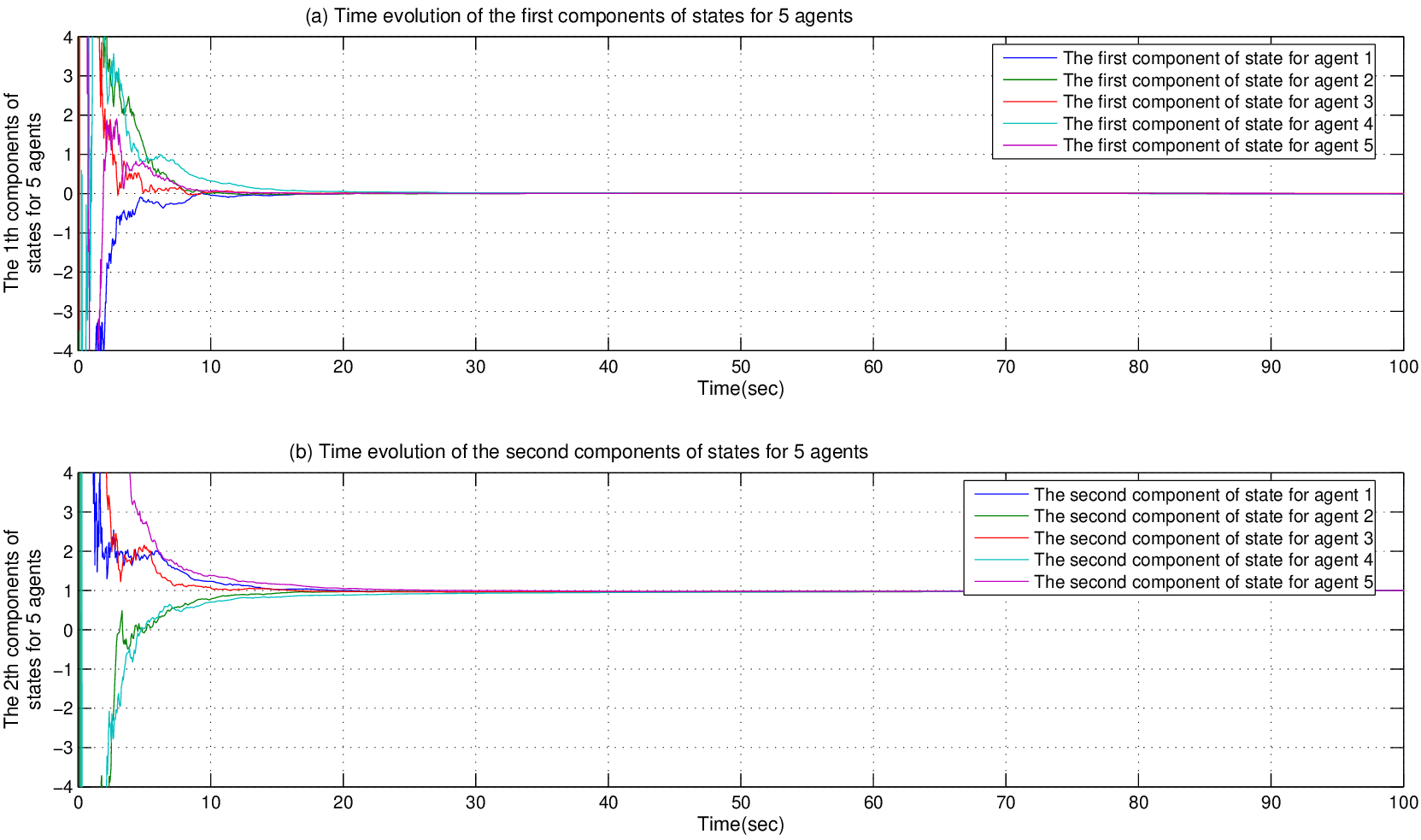}}
\centering{\includegraphics[width=12 cm]{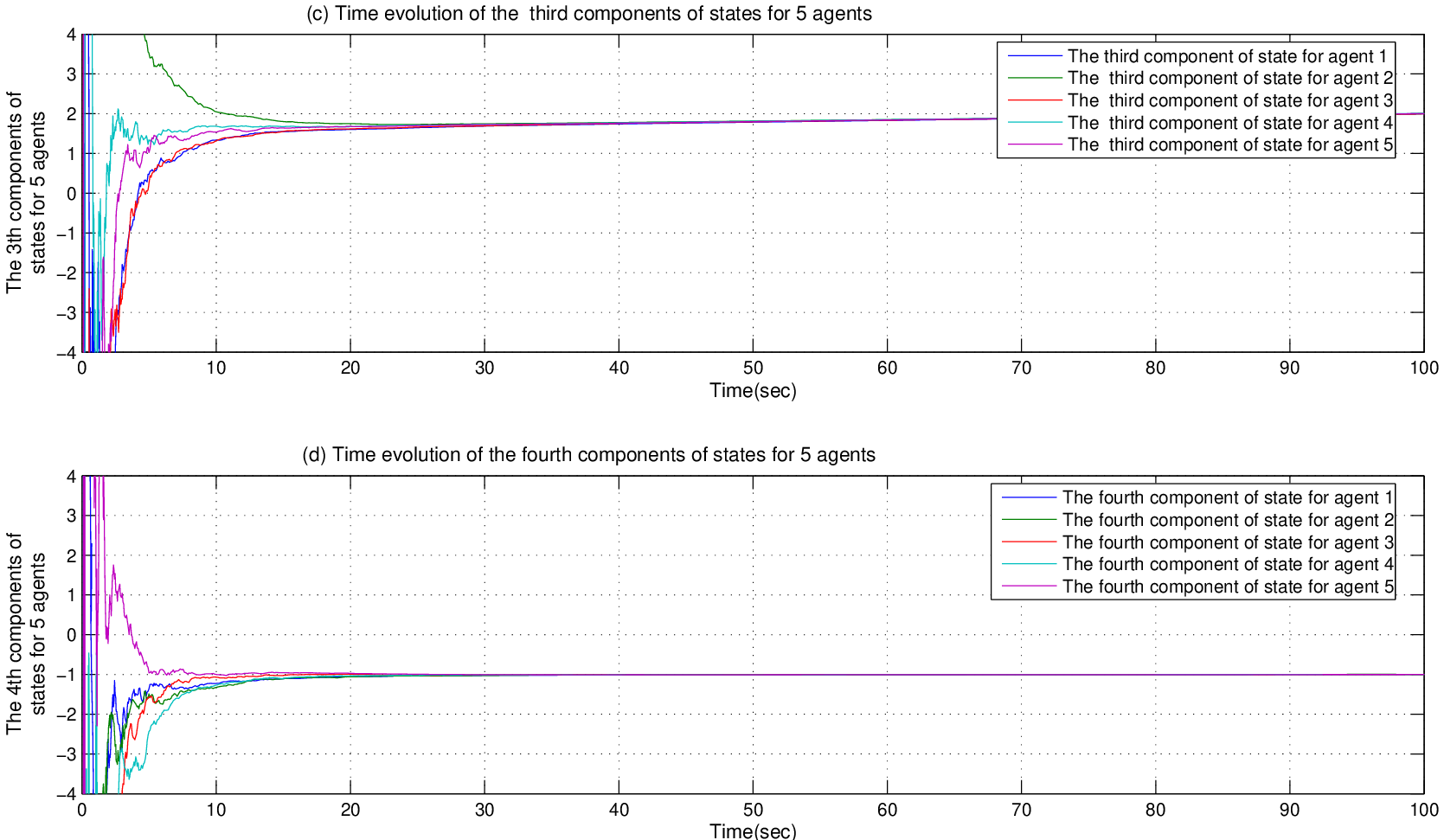}}
\caption{The time evolution of the states for 5  agents.    (a) The first components  of the five agents all converge to $0$; (b)The second components  of the five agents all converge to $1$; (c) The third components  of the five agents all converge to $2$; (d)The fourth components  of the five agents all converge to $-1$.  Therefore, each state of the 5 agents converges to the optimal solution $(0, 1, 2, -1)$.}
\label{EX4}
\end{figure}

We now use our distributed algorithm \eqref{algorithm3} to recover this result.
Note the function corresponding to the equality constraint is not affine which is not satisfied by our theorems. However, since our theorems only give sufficient conditions, this mean that our algorithms may work even if these sufficient conditions are not satisfied. Our simulation shows that this is indeed the case for the above example.
To run the simulation, we  use five agents whose coupling is described by the graph in Fig. \ref{graph}. The state of each agent lies in $\mb R^{4}$. To run simulation under this setup,  define $x=(x_1, x_2, x_3, x_4)^{\top}$ and let $f_1(x)=x_1^2+x_2^2$, $f_2(x)=2x_3^2+x_4^2$, $f_3(x)=-5x_1-5x_2$, $f_4(x)=-21x_3$, $f_5(x)=7x_4$, $g_1(x)=-8+x_1-x_2+x_3-x_4+x_1^2+x_2^2+x_3^2+x_4^2$, $g_2(x)=-10-x_1-x_4+x_1^2+2x_2^2+x_3^2+2x_4^2$, $h_1(x)=-5+2x_1-x_2-x_4+2x_1^2+x_2^2+x_3^2$.
We assume that agent 1 has only access to $f_1, g_1, h_1$, agent 2 has only access to $f_2, g_2$, and agent $i$ has only access to $f_i$ for $i=3, 4, 5$.
The parameters $\eta_{ij}$ in \eqref{algorithm3} is only required to be positive so that we chose them to be $1$ for briefly.   The initial states of the five agents are
$x_1(0)=(3,4,-3,4)^{\top}$, $x_2(0)=(1, -2, 4, 2)^{\top}$,  $x_3(0)=(-3, -4, 3, 3)^{\top}$,  $x_4(0)=(3, 1, 2, -3)^{\top}$,  $x_5(0)=(4, -2, -4, 1)^{\top}$, the initial states for the multipliers corresponding to the inequality constraints of agents 1 and 2 are given as $\l_1(0)=\l_2(0)=3$, the initial state for the multiplier corresponding to the equality constraints of agents 1  is given as$\lb_1(0)=3$, the initial states for the transformed multiplier $\lbb_i$ of the five agents are given as $\lbb_1(0)=\cdots=\lbb_5(0)=(1, 2, 3, 4)^{\top}$. The simulation result, shown in Fig. \ref{EX4}, indicates that the first, second, third, and fourth components of five agents  converge to $0, 1, 2, -1$ respectively. To save place, the time evolutions of the multipliers are not plotted here.

Obviously, we can also use, say 7 agents, to simulate the optimal solution. Which number to use is application depend. Due to space limitation, the simulation result is omitted here.

\section{Conclusions}\label{conclusion}

This paper introduces a projection-free optimization dynamics for nonsmooth composite optimization problems involving both equality and convex inequality constraints. Unlike conventional projection-based approaches, which introduce nonsmooth dynamics and involve complicate stability analysis, the proposed method avoids using projection and ensures a fully smooth optimization process. This smoothness allows for rigorous stability analysis using classical Lyapunov techniques, eliminating the need for nonsmooth analysis tools such as differential inclusions.

Although several studies have also attempted to design constrained optimization algorithms without explicitly relying on projection methods, they instead substitute the projection operator with the relatively easier-to-compute max operator, still rendering  the resulting optimization dynamics nonsmooth in the classical sense.
Our optimization dynamics, which operates without projection, is entirely smooth and facilitates stability analysis through traditional Lyapunov methods, eliminating the necessity for nonsmooth techniques such as differential inclusions. More importantly, unlike the max-based approach which is confined to affine inequality constraints, our method overcomes this limitation and extends to general convex inequality constraints.

Our algorithm also has the ability to achieve exponential convergence and this represents an advancement over existing works that typically ensure only asymptotic or semi-global exponential convergence. On the other hand, our method follows a new line of exponential convergence analysis. Unlike existing methods of proving exponential stability by using the tool of IQC which are valid only when the additional constraints are linear or affine, this paper extends the PAL framework to incorporate general convex inequality constraints and achieves strong convexity-concavity in the primal-dual variables, coming with the capability of dealing with general convex nonlinear constraints.

\section*{ACKNOWLEDGEMENT}
\small
This work is supported by National Natural Science Foundation of China (62363027, 62373180) and Natural Science Foundation of Jiangxi Province (20224BAB202026).

\makesubmdate

%
%
%
%
%
%
%

\makecontacts

\bibliographystyle{plain}
\bibliography{ReferenceOptimization}

\end{document}